\documentclass[reqno,centertags,12pt,a4paper]{amsart}
\usepackage{dsfont}
\usepackage{amscd}
\usepackage{amsmath,tikz}
\usepackage{amssymb}
\usepackage{amsfonts}
\usepackage{hyperref}
\usepackage{enumerate}
\usepackage{latexsym}
\usepackage{cases}
\usepackage{amsthm}
\usepackage{graphicx}
\usepackage{verbatim}
\usepackage{mathrsfs}
\usepackage{color}

\newtheorem{lemma}{Lemma}[section]
\newtheorem{proposition}[lemma]{Proposition}

\newtheorem{remark}[lemma]{Remark}

\newtheorem{theorem}[lemma]{Theorem}
\theoremstyle{definition}
\numberwithin{equation}{section}

\allowdisplaybreaks[2] 

\begin{document}

\title[Schr\"odinger operator along one class of curves]{Sharp pointwise convergence on the Schr\"odinger operator along one class of curves}
\author[Zhenbin~Cao~and~Changxing ~Miao]{Zhenbin~Cao~and~Changxing ~Miao}

\date{\today}

\address{School of Mathematical Sciences, Beijing Normal University, Beijing 100875, China}
\email{11735002@zju.edu.cn}

\address{Institute of Applied Physics and Computational Mathematics, Beijing 100088, China}
\email{miao\_{}changxing@iapcm.ac.cn}


\begin{abstract}
Almost everywhere convergence on the solution of Schr\"odinger equation is an important problem raised by Carleson, which was essentially solved by Du-Guth-Li \cite{DGL} and Du-Zhang \cite{DZ}. In this note, we obtain the sharp pointwise convergence on the Schr\"odinger operator along one class of curves.
\end{abstract}

\maketitle

\section{INTRODUCTION}\label{section1}

We consider the free Schr\"odinger equation:
\begin{align}
 \begin{cases}
  \ i u_{t} - \Delta u = 0, \quad (x,t) \in \mathbb{R}^{n} \times \mathbb{R} , \\
 \  u(x,0)=f(x), \ \ \quad  \,  \  x\in \mathbb{R}^{n}.
 \end{cases}
 \end{align}
Its solution is given by
\begin{equation*}
e^{it\Delta}f(x)= \int_{\mathbb{R}^{n}} e(x \cdot \xi +t|\xi|^{2})\widehat{f}(\xi)d\xi ,
\end{equation*}
where $\widehat{f}$ denotes the Fourier transform of the function $f$, and $e(b):=e^{2\pi ib}$ for each $b \in \mathbb{R}$.

A fundamental problem is determining the  optimal $s$ such that
\begin{equation}\label{aim}
   \lim_{t\rightarrow 0}e^{it\Delta}f(x)=f(x)   \quad \quad  {\rm a.e.}
\end{equation}
for all $f\in H^s(\mathbb{R}^n)$. Carleson \cite{C} first purposed this problem  and  proved \eqref{aim} holds for any $f\in H^{1/4}(\mathbb{R})$. Dahberg and Kenig \cite{DK}
proved the condition $s\ge\frac14$ given by Carleson is sharp. For the situation in higher dimensions, many authors studied this problem \cite{B1,Carbery,Cowling,DG,DGLZ,L,LR,MVV,Sjolin,Vega}.
   In particular, Bourgain \cite{B2} gave counterexamples showing that (\ref{aim}) can fail if $s <\frac{n}{2(n+1)}$.
    Du-Guth-Li \cite{DGL} and  Du-Zhang \cite{DZ} proved \eqref{aim} holds for $s>\frac{n}{2(n+1)}$
     when $n=2$  and $n\ge3$, respectively. Therefore this problem was solved except the endpoint.

     Cho-Lee-Vargas \cite{CLV} considered one class of pointwise convergence problems which is non-tangential convergence to the initial data. By Sobolev embedding, we easily see that non-tangential convergence holds for $s>n/2$. Then Sj\"ogren-Sj\"olin \cite{SS} showed that non-tangential convergence fails for $s\leq n/2$. However, for some special cases of this problem, related results may be improved. Lee-Rogers \cite{LR2} studied the Schr\"odinger operator with the  harmonic oscillator $e^{-it(\Delta+|x|^2)}f(x)$, and showed that after a transformation it can be rewritten as $e^{it\Delta}f(\sqrt{1+t^2}x)$. And they proved, if the curve $\gamma(x,t)$ is $C^1$ function, then the pointwise convergence of $e^{it\Delta}f(\gamma(x,t))$ is essentially equivalent to that of  $e^{it\Delta}f(x)$. Therefore the pointwise convergence of $e^{-it(\Delta+|x|^2)}f(x)$ is essentially equivalent to that of $e^{it\Delta}f(x)$. Their work motivated the study on the relation between the degree of the tangency and regularity when $(x.t)$ approaches to $(x,0)$ tangentially. One model of such problem is
     \begin{equation}\label{ntc}
  \lim_{t\rightarrow0}e^{it\Delta}f(\gamma(x,t))=f(x), \quad  \quad \rm{a.e.}
     \end{equation}
 where $\gamma(x,t)$ is continuous and $\gamma(x,0)=x$. Cho-Lee-Vargas \cite{CLV} considered the following condition: $\gamma(x,t)$ satisfies H\"older condition of order $\alpha$ with $0<\alpha \leq  1$ in $t$:
    \begin{equation}\label{con1}
    |\gamma(x,t)-\gamma(x,t')|\leq C|t-t'|^\alpha,
    \end{equation}
     and bi-Lipschitz in $x$:
     \begin{equation}\label{con2}
     	C_1|x-y|\leq |\gamma(x,t)-\gamma(y,t)|\leq C_2|x-y|.
     \end{equation}
When $n=1$, they proved (\ref{ntc}) holds for $s>\max(1/2-\alpha,1/4)$ if $0<\alpha \leq  1$ by the stationary phase method, and this result is sharp up to the endpoint. Ding-Niu \cite{DN} studied weighted $L^p$ maximal estimates on $e^{it\Delta}f(\gamma(x,t))$. As a corollary, they showed that (\ref{ntc}) holds for $s\geq 1/4$ if $1/2\leq \alpha \leq  1$. When $n=2$, Li-Wang \cite{LW} proved (\ref{ntc}) holds for $s>3/8$ if $1/2 \leq \alpha <1$ by the polynomial partitioning method.

 A simple example satisfying (\ref{con1}) and (\ref{con2}) is $\gamma(x,t)=x+t^\alpha \mu$ for some $\mu \in \mathbb{R}^n \backslash \{0\}$. However, only studying this case is not enough to complete our induction. Thus we consider one class of curves $\gamma(x,t)$, which defined as
     \begin{equation}\label{con3}
\gamma(x,t):=x+ t^\alpha \mu(t),  \quad \quad  t\geq 0.
     \end{equation}
     Here  $\mu(t) \in C^\infty([0,\infty) )\times ...\times C^\infty([0,\infty) ),$ and
     \begin{equation}\label{eab}
     	|\mu^{(k)}(t)|\lesssim_k t^{-k}, \quad \quad  \forall ~k \in \mathbb{N}.
     \end{equation}
We denote the collection of the curves $\gamma(x,t)$ satisfying (\ref{con3}) by $\mathcal{D}$. Define the operator associated with such $\gamma(x,t)$ as
     $$  T_\alpha^\gamma f(x,t):=  e^{it\Delta}f(\gamma(x,t))=\int_{\mathbb{R}^n} e(x\cdot \xi+t^\alpha\mu(t) \cdot \xi +t|\xi|^2)\widehat{f}(\xi)d\xi. $$
     Then  main result is the following:

\begin{theorem}\label{th1}
Let $n\geq 2$, $1/2 \leq \alpha <1$ and $\gamma(x,t) \in \mathcal{D}$. For every $f\in H^s(\mathbb{R}^n)$ with $s>\frac{n}{2(n+1)}$,
	$$   \lim_{t\rightarrow 0+}  T_\alpha^\gamma f(x,t)=f(x), \quad\quad  \rm{a.e.} $$
	and the range of $s$ is sharp up to the endpoint.
\end{theorem}

The sharp property in Theorem \ref{th1} can be checked through Bourgain's counterexample in \cite{B2}. We refer to \cite[Theorem 1.12]{LW}. On the positive part of Theorem \ref{th1}, we can reduce it to the following estimate.

\begin{theorem}\label{th2}
Let $n\geq 2$, $1/2 \leq \alpha <1$ and $\gamma(x,t) \in \mathcal{D}$. For every $f\in H^s(\mathbb{R}^n)$ with $s>\frac{n}{2(n+1)}$, there exists a constant $C_s$ such that
\begin{equation}\label{th2 eq}
\left\|  \sup_{0<t\leq 1} |T_\alpha^\gamma f| \right\|_{L^2(B^n(0,1))}   \leq C_s \|f\|_{H^s(\mathbb{R}^n)}.
\end{equation}
\end{theorem}

Via Littlewood-Paley decomposition, if $f$ is Fourier supported on $A(R)=\{   \xi \in \mathbb{R}^n : |\xi|\sim R \}$, (\ref{th2 eq}) boils down to the bound
	$$   \left\|  \sup_{0<t\leq 1} |T_\alpha^\gamma f| \right\|_{L^2(B^n(0,1))}   \leq C_\epsilon R^{\frac{n}{2(n+1)}+\epsilon}\|f\|_{L^2}. $$
After time localization lemma \cite[Lemma 2.1]{CLV}, this bound is further reduced to
	$$   \left\|  \sup_{\lambda<t\leq \lambda+R^{-1}} |T_\alpha^\gamma f| \right\|_{L^2(B^n(0,1))}   \leq C_\epsilon R^{\frac{n}{2(n+1)}+\epsilon}\|f\|_{L^2},  \quad \quad \forall~ \lambda \in \frac{1}{R}\mathbb{Z} \cap [0,1-\frac{1}{R}]. $$
For $\gamma(x,t) \in \mathcal{D}$, we define
$$ \gamma_R(x,t):=x+R^{1-2\alpha}t^\alpha \mu\Big(\frac{t}{R^2}\Big).  $$
Then using parabolic rescaling, it suffices to show the following estimate.
\begin{theorem}\label{th3}
Let $n\geq 2$, $1/2 \leq \alpha <1$ and $\gamma(x,t) \in \mathcal{D}$.
For any $\epsilon>0$, there exists a constant $C_\epsilon$ such that
\begin{equation}\label{th3 eq}
\left\|  \sup_{\lambda<t\leq \lambda+R} |T_\alpha^{\gamma_R} f| \right\|_{L^2(B^n(0,R))}   \leq C_\epsilon R^{\frac{n}{2(n+1)}+\epsilon}\|f\|_{L^2}
\end{equation}
holds for all $R\geq 1$, every $\lambda \in R\mathbb{Z} \cap [0,R^2-R]$, and all $f$ with supp$\widehat{f} \subset B^n(0,1).$
\end{theorem}

Finally, we can reduce Theorem \ref{th3} to the following result which says that only considering $\lambda=0$ is enough.

\begin{theorem}\label{th4}
	Let $n\geq 2$, $1/2 \leq \alpha <1$ and $\gamma(x,t) \in \mathcal{D}$.
	For any $\epsilon>0$, there exists a constant $C_\epsilon$ such that
	\begin{equation}\label{th4 eq}
		\left\|  \sup_{0<t\leq R} |T_\alpha^{\gamma_R} f| \right\|_{L^2(B^n(0,R))}   \leq C_\epsilon R^{\frac{n}{2(n+1)}+\epsilon}\|f\|_{L^2}
	\end{equation}
	holds for all $R\geq 1$ and all $f$ with supp$\widehat{f} \subset B^n(0,1).$
\end{theorem}

For the case $e^{it\Delta}f(x)$, the feature of Theorem \ref{th4} implying Theorem \ref{th3} is obvious since $e^{it\Delta}f(x)$ is translation invariant on the variable $t$. As for $e^{it\Delta}f(\gamma(x,t))$, where $\gamma(x,t)=x+t^\alpha \mu$ for some $\mu \in \mathbb{R}^n \backslash \{0\}$, this property is invalid. To overcome this difficulty, we consider one larger class of curves $\gamma(x,t) \in \mathcal{D}$, which maintains this property. We leave the proof of Theorem \ref{th4} implying Theorem \ref{th3} to Section \ref{sec2}. Another advantage of restricting $\gamma(x,t) \in \mathcal{D}$ is that $\gamma(x,t)$ is smooth for $t$ away from 0. This feature allows us to use the decoupling inequality of variable coefficient version, which demands the phase function is $C^2$ at least. The details of this part will be spread out in Section \ref{sec3}.

\vskip0.2cm

\noindent\textbf{Notation.} If $X$ is a finite set, we use $\#X$ to denote its cardinality. If $X$ is a measurable set,
we use $|X|$ to denote its Lebesgue measure. $C_{\epsilon}$  denotes a constant which depends on
 $\epsilon$.  Write $A\lesssim B$ or $A=O(B)$ to
mean that there exists a constant $C$ such that $A\leq CB$. We use $B^{n+1}(c,r)$ to denote a ball centered at $c$ with radius $r$ in $\mathbb{R}^{n+1}$. We abbreviate $B^{n+1}(c,r)$ to $B(c,r)$ if $\mathbb{R}^{n+1}$ is clear in the context.

\section{Preparation and basic induction}\label{sec2}

Main ingredients of our proof include locally constant property, dyadic pigeonholing, broad-narrow analysis, parabolic rescaling and induction on scale, which has same techniques as in  \cite{DGL,DGLZ,DZ}. Among them, locally constant property is a basic tool. It can reduce mixed norm in Theorem \ref{th4} to uniform norm on a sparse set. Locally constant property says that if a function $f$ is Fourier supported on a $\rho^{-1}$-ball, then we can view $|f|$ essentially as constant on every ball at scale $\rho$. There are several versions, such as Bernstein's inequality, stability lemma \cite[Lemma 2.1]{T2}, \cite[Lemma 6.1]{GWZ}, and corresponding variable coefficient version \cite[Lemma 5.8]{GHI}. Unfortunately, these results can't be used on $T_\alpha^{\gamma_R} f$ (and $T_\alpha^{\gamma} f$), since it has no compact Fourier support, and its phase function is not smooth on the variable $t$. However, compact Fourier support is not the necessary condition to make locally constant property take effect. Another explanation of this property is that the term $e(a(x,y))$ can be removed from the integration $\int e(a(x,y)) g(y)dy$ if $|a(x,y)| \lesssim O(1)$. We give the following proposition.

\begin{proposition}\label{lcp}
Let $d\in \mathbb{N}^+$, $M,N>0$, $x_0,y_0 \in \mathbb{R}^d$ and $q \geq 1$. Define
$$  Sg(x)=   \int e(a(x,y)) K(x,y) g(y)dy$$
 and
$$   \tilde{S}g(x)=\int e(a(x_0,y))K(x,y) g(y)dy,$$
where $x \in B^d(x_0,M),$ supp$g \subset B^d(y_0,N)$. Let $\tilde{a}=a(x,y)-a(x_0,y)$.

\noindent\rm{(a)} If $\tilde{a}(x,y)$ is smooth on $x$, and
\begin{equation}\label{lcp c1}
		\sup_{\substack{     x \in B(x_0,M)   \\    y\in  B(y_0,N)    }}|\partial_x^\beta \tilde{a}(x,y)| \lesssim_\beta M^{-|\beta|}, \quad \quad \beta \in \mathbb{N}^d \text{~with~} 1 \leq |\beta|\leq d+1,
\end{equation}
then
\begin{equation}\label{lcp c3}
	\|Sg\|_{L^q(B(x_0,M))} \lesssim \sum_{k\in \mathbb{Z}^d} (1+|k|)^{-(d+1)} \left\|\tilde{S}\left( c_k g   \right)\right\|_{L^q(B(x_0,M))}
\end{equation}
for some choices of smooth functions $c_k$ satisfying the uniform bound $\|c_k\|_{L^\infty(B(y_0,N))} \lesssim 1.$

\noindent\rm{(b)} If $\tilde{a}(x,y)$ is smooth on $y$, and
\begin{equation}\label{lcp c2}
	\sup_{\substack{     x \in B(x_0,M)   \\    y\in  B(y_0,N)    }}|\partial_y^\beta \tilde{a}(x,y)| \lesssim_\beta N^{-|\beta|}, \quad \quad \beta \in \mathbb{N}^d \text{~with~} 1 \leq |\beta|\leq d+1,
\end{equation}
then
\begin{equation}\label{lcp c4}
	\|Sg\|_{L^q(B(x_0,M))} \lesssim \sum_{k\in \mathbb{Z}^d} (1+|k|)^{-(d+1)} \left\|\tilde{S}\Big(g  e\Big(  \frac{ky}{N}\Big)  \Big)\right\|_{L^q(B(x_0,M))}.
\end{equation}
\end{proposition}

\noindent \textit{Proof.} We only prove (b) since the arguments of two cases are parallel. Take a bump function $\varphi$ supported on $(-\pi,\pi)^d$, and $\varphi=1$ on $B^d(0,1)$. Then
$$   Sg(x)=\int \varphi\left(\frac{y-y_0}{N}\right)e(\tilde{a}(x,y))e(a(x_0,y))K(x,y) g(y)dy . $$
By expanding $\varphi(\frac{y-y_0}{N})e(\tilde{a}(x,y))$ as a Fourier series in the variable $y$, one has
\begin{align*}
	\varphi\left(\frac{y-y_0}{N}\right)e(\tilde{a}(x,y))=\sum_{k\in \mathbb{Z}^d} c_k(x) e\left( \frac{k(y-y_0)}{N}\right),
\end{align*}
where each $c_k(x)$ satisfies
$$\sup_{x\in B(x_0,M)}|c_k(x)|\lesssim  (1+|k|)^{-(d+1)}      $$
due to (\ref{lcp c2}).
Then
$$   Sg(x)=   \sum_{k\in \mathbb{Z}^d} c_k(x) e(-ky_0/N)\tilde{S}(g e(ky/N))(x),$$
and (\ref{lcp c4}) follows immediately by Minkowski's inequality.

\qed

 (\ref{lcp c3}) and (\ref{lcp c4}) imply that the estimate of $Sg$ can be reduce to that of $\tilde{S}g$. From now on, for convenient, we write (\ref{lcp c3}) and (\ref{lcp c4}) as
$$   |Sg| \sim |\tilde{S}g|.  $$
\vskip0.5cm

Now we show how Theorem \ref{th3} follows from Theorem \ref{th4}.
If $\lambda=0$, (\ref{th3 eq}) and (\ref{th4 eq}) coincide. If $R \leq \lambda<R^{1+\epsilon}$, set $w=2R^\epsilon$, then
$$   \left\|  \sup_{\lambda<t\leq \lambda+R} |T_\alpha^{\gamma_R} f(x,t)| \right\|_{L^2(B_R)} \\
\leq    \left\|  \sup_{0<t\leq wR} |T_\alpha^{\gamma_R} f(x,t)| \right\|_{L^2(B_{R})} . $$
After parabolic rescaling $ t\rightarrow t/w^2, x\rightarrow x/w$, we have
$$ \left\|  \sup_{0<t\leq wR} |T_\alpha^{\gamma_R} f(x,t)| \right\|_{L^2(B_{R})}= R^{O(\epsilon)}\left\|  \sup_{0<t\leq R/w} |T_\alpha^{\gamma_{R/w}} g(x,t)| \right\|_{L^2(B_{R/w})}    $$
for some $g$ satisfying $\|g\|_{L^2}=\|f\|_{L^2}$ and supp$\widehat{g} \subset B(0,w) $. We decompose $B(0,w)$ into unit balls $B(v_j,1)$. Then write $g=\sum_j g_j$, where each $g_j$ is Fourier supported on $B(v_j,1)$. Note  $\# j\sim w^n$ and $|v_j|\leq w$. Define $h_j$ by
$$  \widehat{h_j}(\xi):=\widehat{g_j}(\xi+v_j).   $$
Thus by Theorem \ref{th4} we obtain
\begin{align*}
	&	\left\|  \sup_{0<t\leq R/w} |T_\alpha^{\gamma_{R/w}} g(x,t)| \right\|_{L^2(B_{R/w})}  \\
	\leq & \sum_j  	\left\|  \sup_{0<t\leq R/w} |T_\alpha^{\gamma_{R/w}} g_j(x,t)| \right\|_{L^2(B_{R/w})}  \\
	= & \sum_j  	\left\|  \sup_{0<t\leq R/w} |T_\alpha^{\gamma_{R/w}} h_j(x+2tv_j,t)| \right\|_{L^2(B_{R/w})}  \\
	\leq & \sum_j  	\left\|  \sup_{0<t\leq R/w} |T_\alpha^{\gamma_{R/w}} h_j(x,t)| \right\|_{L^2(B_{3R})}   \\
	\leq  &  \sum_{j,k}  	\left\|  \sup_{0<t\leq R/w} |T_\alpha^{\gamma_{R/w}} h_j(x,t)| \right\|_{L^2(B(x_k,R/w))}  \\
	= &   \sum_{j,k}  	\left\|  \sup_{0<t\leq R/w} |T_\alpha^{\gamma_{R/w}} h_{j,k}(x,t)| \right\|_{L^2(B(0,R/w))} \\
	\lesssim & R^{\frac{n}{2(n+1)}+O(\epsilon)} \sum_{j,k}   \|h_{j,k}\|_{L^2}  \\
	\lesssim & R^{\frac{n}{2(n+1)}+O(\epsilon)} \|f\|_{L^2}.
\end{align*}
Here in the fifth line we break $B_{3R}$ into finitely overlapping balls of the form $B(x_k,R/w)$, and $h_{j,k}$ is given by
$$  \widehat{h_{j,k}}(\xi):=e(x_k \cdot \xi) \widehat{h_j}(\xi).  $$
Combining all estimates above, we complete the argument of this case.

If $\lambda\geq R^{1+\epsilon}$, after translation,
$$ \left\|  \sup_{\lambda<t\leq \lambda+R} |T_\alpha^{\gamma_R} f(x,t)| \right\|_{L^2(B_R)}= \left\|  \sup_{0<t\leq R} |T_\alpha^{\gamma_R} f(x,t+\lambda)| \right\|_{L^2(B_R)}.$$
Note
\begin{align*}
	&	|T_\alpha^{\gamma_R} f(x,t+\lambda)| \\
	= &\left|\int e\left[  x\cdot \xi+R^{1-2\alpha}(t+\lambda)^\alpha\mu\left(\frac{t+\lambda}{R^2}\right)\cdot\xi+t|\xi|^2  \right]\left[ e(\lambda|\xi|^2)\widehat{f}(\xi)   \right]d\xi\right|   .
\end{align*}
For $0<t \leq R$ and $\lambda \geq R^{1+\epsilon}$, it follows from Taylor's formula that
\begin{align*}
&(t+\lambda)^\alpha\mu\left(\frac{t+\lambda}{R^2}\right) \\
 =& \lambda^\alpha    \left[    \sum_{i=0}^k c_i\left(\frac{t}{\lambda}\right)^i+R_{k+1}(\lambda,t)\right]\mu\left(\frac{t+\lambda}{R^2}\right)     \\
 =&  t^\alpha    \left[    \sum_{i=1}^k c_i\left(\frac{t}{\lambda}\right)^{i-\alpha}\mu\left(\frac{t+\lambda}{R^2}\right)\right]  + \lambda^\alpha \mu\left(\frac{t+\lambda}{R^2}\right)  +\lambda^\alpha R_{k+1}(\lambda,t)\mu\left(\frac{t+\lambda}{R^2}\right)  \\
 =&  t^\alpha    \left[    \sum_{i=1}^k c_i\left(\frac{t}{\lambda}\right)^{i-\alpha}\mu\left(\frac{t+\lambda}{R^2}\right)\right]  + \lambda^\alpha \left[\sum_{j=0}^k d_j \mu^{(j)}\left(\frac{\lambda}{R^2}\right) \left(\frac{t}{R^2}\right)^j +R'_{k+1}(\lambda,t)\right]  \\
 & \quad\quad\quad\quad\quad\quad\quad\quad\quad\quad\quad\quad\quad\quad\quad\quad\quad\quad\quad\quad +\lambda^\alpha R_{k+1}(\lambda,t)\mu\left(\frac{t+\lambda}{R^2}\right)  \\
 =& t^\alpha    \left[    \sum_{i=1}^k c_i\left(\frac{t}{\lambda}\right)^{i-\alpha}\mu\left(\frac{t+\lambda}{R^2}\right)     + \sum_{j=1}^k d_j \mu^{(j)}\left(\frac{\lambda}{R^2}\right) \left(\frac{t}{R^2}\right)^j   \left(\frac{\lambda}{t}\right)^\alpha       \right]  +  \lambda^\alpha\mu\left(\frac{\lambda}{R^2}\right)   \\
 &  \quad\quad\quad\quad\quad\quad\quad\quad\quad\quad\quad\quad   +\lambda^\alpha R_{k+1}(\lambda,t)\mu\left(\frac{t+\lambda}{R^2}\right) +   \lambda^\alpha R'_{k+1}(\lambda,t) \\
 :=& t^\alpha \tilde{\mu}\left(\frac{t}{R^2}\right) + \lambda^\alpha\mu\left(\frac{\lambda}{R^2}\right)  +\lambda^\alpha R_{k+1}(\lambda,t)\mu\left(\frac{t+\lambda}{R^2}\right) +   \lambda^\alpha R'_{k+1}(\lambda,t) ,
\end{align*}
where
$$   \tilde{\mu}(\tilde{t}):= \sum_{i=1}^k c_j\left(\frac{\tilde{t}}{\tilde{\lambda}}\right)^{i-\alpha}\mu(\tilde{t}+\tilde{\lambda}) +\sum_{j=1}^k d_j \mu^{(j)}(\tilde{\lambda})\tilde{\lambda}^\alpha \tilde{t}^{j-\alpha}, \quad \tilde{\lambda}=\frac{\lambda}{R^2},  $$
$R_{k+1}$ and $R'_{k+1}$ denote Taylor remainders. Then $	|T_\alpha^{\gamma_R} f(x,t+\lambda)| $ can be rewritten as
\begin{align*}
&\bigg|\int e\left[  x\cdot \xi+R^{1-2\alpha}\left( t^\alpha \tilde{\mu}\left(\frac{t}{R^2}\right) +\lambda^\alpha R_{k+1}(\lambda,t)\mu\left(\frac{t+\lambda}{R^2}\right) +   \lambda^\alpha R'_{k+1}(\lambda,t)\right)\cdot\xi+t|\xi|^2  \right]   \\
&  \quad \quad\quad \quad \quad \quad \quad \quad \quad \quad \quad \quad \quad \quad  \left[ e\left(R^{1-2\alpha}\lambda^\alpha\mu\left(\frac{\lambda}{R^2}\right)\cdot\xi+\lambda|\xi|^2\right)\widehat{f}(\xi)   \right]d\xi\bigg|.
\end{align*}
By a direct calculation, we get $\tilde{\mu}$ satisfies (\ref{eab}). On the other hand, by taking $k\in \mathbb{N}$ with $k>(1-\alpha)\frac{1-\epsilon}{\epsilon}$, we can remove the terms  $R^{1-2\alpha} \lambda^\alpha R_{k+1}(\lambda,t)\mu(\frac{t+\lambda}{R^2})\cdot\xi$ and $R^{1-2\alpha}\lambda^\alpha R'_{k+1}(\lambda,t)\cdot\xi$ from above phase function by Proposition \ref{lcp}. Finally, we obtain
\begin{align*}
&	|T_\alpha^{\gamma_R} f(x,t+\lambda)|   \\
	\sim &\left|\int e\left[x\cdot\xi  +R^{1-2\alpha}t^\alpha\tilde{\mu}\left(\frac{t}{R^2}\right)\cdot\xi  +t|\xi|^2\right]\left[ e\left(R^{1-2\alpha}\lambda^\alpha\mu\left(\frac{\lambda}{R^2}\right)\cdot\xi+\lambda|\xi|^2\right)\widehat{f}(\xi)   \right]d\xi\right| \\
	:=& |T_\alpha^{\gamma_R} f_\lambda(x,t)| ,
\end{align*}
where
$$\widehat{f}_\lambda(\xi):=e\left(R^{1-2\alpha}\lambda^\alpha\mu\left(\frac{\lambda}{R^2}\right)\cdot\xi+\lambda|\xi|^2\right)\widehat{f}(\xi) .$$
We use Theorem \ref{th4} to $T_\alpha^{\gamma_R} f_\lambda $, then the result follows.

\section{The proof of Theorem \ref{th4}}\label{sec3}

In this section, we start to prove Theorem \ref{th4}. We can view $|T_\alpha^{\gamma_R} f|$ as constant on each unit ball by using Proposition \ref{lcp}. Then using dyadic pigeonholing, (\ref{th4 eq}) can be reduced to
	\begin{equation}\label{aaaaaa}
	\left\| T_\alpha^{\gamma_R} f \right\|_{L^2(X)}   \lesssim R^{\frac{n}{2(n+1)}+\epsilon}\|f\|_{L^2}.
\end{equation}
Here $X$ denotes a union of unit balls in $B^n(0,R) \times [0.R]$ satisfying the property that each vertical thin tube of dimensions $1 \times...\times 1\times R$ contains exactly one unit ball in $X$. Du-Zhang \cite{DZ} produced one inductive mechanism called the fractal $L^2$ restriction estimate to deal with (\ref{aaaaaa}). More precisely, they introduced a parameter $\gamma$, motivated by Falconer distance set problem, which measures the sparse property on the unit balls in $X$. Then the sparse property of $\gamma$, together with square root cancellation property of the $l^2$ decoupling inequality, makes their proof work. Via the similar argument as them, we can reduce (\ref{aaaaaa}) to the following proposition. Denote $B^\ast_R=B^n(0,R) \times [0,R]$.

\begin{proposition}\label{pro}
Let $n \geq 2$, $1/2 \leq \alpha <1$ and $\gamma(x,t) \in \mathcal{D}$. For any $0<\epsilon <1/100$, there exist constant $C_\epsilon$ and $\delta=\min(\epsilon^{100},\frac{1-\alpha}{100})$ such that the following holds for all $R \geq 1$ and all $f$ with supp$\widehat{f} \subset B^n(0,1)$. Let $p=\frac{2(n+1)}{n-1}$. Suppose that $Y=\cup_{k=1}^M B_k$ is a union of lattice $K^2$-cubes in $B^\ast_R$, where $K=R^\delta$. Suppose that
	$$    \|T_\alpha^{\gamma_R} f\|_{L^p(B_k)} \text{~is~essentially~a~dyadic~number~in~}   k=1,2,...,M.   $$
Let $1\leq \lambda\leq n+1$ and $\nu$ be defined by
	$$  \nu:=\max_{\substack{  B^{n+1}(x',r)\subset B^\ast_R \\ x'\in \mathbb{R}^{n+1},r\geq K^2 }} \frac{\#\{ B_k:B_k \subset B(x',r)\}}{r^\lambda}.  $$
Then
	$$   \left\|   T_\alpha^{\gamma_R} f \right\|_{L^p(Y)}   \leq C_\epsilon M^{-\frac{1}{n+1}}\nu^{\frac{1}{n+1}}R^{\frac{\lambda}{2(n+1)}+\epsilon}\|f\|_{L^2}. $$
\end{proposition}

\begin{remark}
Du-Zhang \cite{DZ} introduced another parameter to compute the number of unit cubes in $Y$ that intersect a given lattice $R^{1/2}$ cube. This more sophisticated result needs the refined Strichartz estimate on $e^{it\Delta}f(x)$ produced by Du-Guth-Li-Zhang \cite{DGLZ}. We can also build the refined Strichartz estimate on $e^{it\Delta}f(\gamma(x,t))$ by following the argument from \cite{DGLZ}. However, as they say in \cite{DZ}, this parameter is not necessary to obtain the relevant pointwise convergence result. For this reason, we prove  above weaker version without this parameter.
\end{remark}

\noindent \textit{Proof.} We decompose $B^n(0,1)$ in the frequency space into $K^{-1}$-cubes $\tau$. Then write $f=\sum_{\tau} f_\tau$, where $\widehat{f_\tau}=\widehat{f}|_\tau$. For a $K^2$-cube $B$ in $Y$, we define associated significant set as
$$    \mathcal{S}(B) :=\Big\{ \tau: \| T_\alpha^{\gamma_R} f _\tau\|_{L^p(B)} \geq \frac{1}{100\# \tau} \|  T_\alpha^{\gamma_R} f  \|_{L^p(B)} \Big\}.  $$
A basic property is
$$  \Big\|\sum_{\tau \in  \mathcal{S}(B)} T_\alpha^{\gamma_R} f _\tau\Big\|_{L^p(B)} \sim \| T_\alpha^{\gamma_R} f \|_{L^p(B)}. $$
We say $B\subset Y$ is broad if there exist $\tau_1,...,\tau_{n+1} \in  \mathcal{S}(B)$ such that for any $v_j \in G(\tau_j)$,
\begin{equation}\label{broad condition}
	|v_1 \wedge v_2 \wedge...\wedge v_{n+1}  |\gtrsim K^{-n},
\end{equation}
where $G(\tau)$ is defined by
$$    G(\tau):=\left\{ \frac{(-2\xi,1)}{|(-2\xi,1)|}\in S^n:\xi \in \tau   \right\}. $$
Otherwise, we say $B$ is narrow. We denote the union of broad $B \subset Y$ by $Y_{\text{broad}}$ and the union of narrow $B \subset Y$ by $Y_{\text{narrow}}$. We call it the broad case if $Y_{\text{broad}}$ contains $\geq M/2$ many $K^2$-cubes, and the narrow case otherwise.

\vskip0.5cm

\noindent \textbf{Broad case.} For each broad $B$, by the definition of $\mathcal{S}(B)$, there exist $\tau_1,...,\tau_{n+1} \in  \mathcal{S}(B)$  satisfying (\ref{broad condition}), and
\begin{equation}\label{broad eq1}
	\|T_\alpha^{\gamma_R} f \|^p_{L^p(B)} \leq K^{O(1)} \prod_{j=1}^{n+1}\left( \int_B |T_\alpha^{\gamma_R} f _{\tau_j}|^p \right)^{\frac{1}{n+1}}.
\end{equation}
Set $B=B(x_B, K^2)$. For applying locally constant property, we decompose $B$ to balls of the form $B(x_B +v,2)$, where $v \in B(0,K^2)\cap \mathbb{Z}^{n+1}$. We choose $v_j \in B(0,K^2) \cap \mathbb{Z}^{n+1}$ such that $\|T_\alpha^{\gamma_R} f _{\tau_j}\|_{L^\infty(B)}$ is attained in $B(x_B+v_j,2)$. Define $v_j=(x_j,t_j)$ and
$$  \widehat{f_{\tau_j,v_j}}(\xi):= \widehat{f_{\tau_j}}(\xi)  e(x_j\cdot\xi+t_j|\xi|^2). $$
Thus
$$  \int_B |T_\alpha^{\gamma_R} f _{\tau_j}(x,t)|^p \leq K^{O(1)}\int_{B(x_B+v_j,2)} |T_\alpha^{\gamma_R} f _{\tau_j}(x,t)|^p = K^{O(1)}\int_{B(x_B,2)} |T_\alpha^{\gamma_R} f _{\tau_j}(x+x_j,t+t_j)|^p. $$
Note
\begin{align*}
&	|T_\alpha^{\gamma_R} f _{\tau_j}(x+x_j,t+t_j)|    \\
	=& \left|\int e\left[(x+x_j)\cdot \xi+R^{1-2\alpha}(t+t_j)^\alpha\mu\left(  \frac{t+t_j}{R^2}\right)\cdot \xi+(t+t_j)|\xi|^2\right]              \widehat{f_{\tau_j}}(\xi) d\xi  \right| \\
	=&  \left|\int e\left[x\cdot \xi+R^{1-2\alpha}(t+t_j)^\alpha\mu\left(  \frac{t+t_j}{R^2}\right)\cdot \xi+t|\xi|^2\right]              \widehat{f_{\tau_j,v_j}}(\xi) d\xi  \right| \\
	=& \bigg|  \int e\left[ R^{1-2\alpha} \left(  (t+t_j)^\alpha\mu\left( \frac{t+t_j}{R^2}\right)-t^\alpha\mu\left(\frac{t}{R^2}\right) \right)\cdot (\xi-c_{\tau_j})    \right]  \\
	& \quad \quad  \quad \quad  \quad \quad  \quad \quad  \cdot e\left[x\cdot \xi+R^{1-2\alpha}t^\alpha\mu\left(  \frac{t}{R^2}\right)\cdot \xi+t|\xi|^2\right]              \widehat{f_{\tau_j,v_j}}(\xi) d\xi  \bigg|  \\
	\sim & |T_\alpha^{\gamma_R} f _{\tau_j,v_j}(x,t) |,
\end{align*}
where $c_{\tau_j}$ denotes the center of $\tau_j$. Here we use Proposition \ref{lcp} to remove the first term in the above phase function. Therefore,
\begin{align*}
	\|T_\alpha^{\gamma_R} f \|^p_{L^p(B)} &  \leq K^{O(1)} \prod_{j=1}^{n+1} \left( \int_{B(x_B,2)} |T_\alpha^{\gamma_R} f _{\tau_j,v_j}|^p  \right)^{\frac{1}{n+1}}   \\
	& \sim K^{O(1)} \int_{B(x_B,2)} \prod_{j=1}^{n+1} |T_\alpha^{\gamma_R} f _{\tau_j,v_j}|^{\frac{p}{n+1}} .
\end{align*}
Since there are only $K^{O(1)}$ choices for $\tau_j,v_j$, by pigeonholing, there exist $\tilde{\tau_j},\tilde{v_j}$ such that above inequality holds for at least $\geq K^{-C}M$ broad $B$. Now fix $\tilde{\tau_j},\tilde{v_j}$.  We abbreviate $f_{\tilde{\tau_j},\tilde{v_j}}$ to $f_j$, and denote the collection of remaining broad $B$ by $\mathcal{B}$. Next we sort $B \in \mathcal{B}$ by the value $ \|\prod_{j=1}^{n+1} |T_\alpha^{\gamma_R} f _j|^{\frac{1}{n+1}}\|_{L^\infty(B(x_B,2))}$: for dyadic number $A$, define
$$ \mathbb{Y}_A:=\left\{B \in \mathcal{B}: \left\|\prod_{j=1}^{n+1} |T_\alpha^{\gamma_R} f _j|^{\frac{1}{n+1}}\right\|_{L^\infty(B(x_B,2))} \sim A\right\} .  $$
Let $Y_{A}$ be the union of the $K^2$-cubes $B$ in $\mathbb{Y}_{A}$. Without loss of generality, assume that $\|f\|_{L^2}=1$. We can further assume that $R^{-C} \leq A \leq 1$ for some constant $C$. So there are only $O(\log R) \leq O(K)$ choices on $A$. By dyadic pigeonholing, there exists a constant $\tilde{A}$ such that
$$  \# \{  B: B\subset Y_{\tilde{A}}  \} \gtrsim K^{-C} \#  \mathcal{B}.  $$
Now we fix $\tilde{A}$, and denote $Y_{\tilde{A}}$ by $Y'$. Combining all above estimates, we obtain
\begin{align*}
	\|T_\alpha^{\gamma_R} f \|^p_{L^p(Y_{\text{broad}})} & \leq 	K^{O(1)}\|T_\alpha^{\gamma_R} f \|^p_{L^p(Y')}   \\
	& \leq K^{O(1)} \left\| \prod_{j=1}^{n+1} |T_\alpha^{\gamma_R} f_j |^{\frac{1}{n+1}}  \right\|_{L^p(\cup_{B\subset Y'} B(x_B,2))} \\
	& \sim K^{O(1)} M^{-\frac{1}{2(n+1)}} \left\| \prod_{j=1}^{n+1} |T_\alpha^{\gamma_R} f _j|^{\frac{1}{n+1}}  \right\|_{L^{\frac{2(n+1)}{n}}(\cup_{B\subset Y'} B(x_B,2))}    \\
	& \leq   K^{O(1)} M^{-\frac{1}{2(n+1)}}  \left\| \prod_{j=1}^{n+1} |T_\alpha^{\gamma_R} f _j|^{\frac{1}{n+1}}  \right\|_{L^{\frac{2(n+1)}{n}}(B_R^\ast)} \\
	& \lesssim    K^{O(1)} M^{-\frac{1}{2(n+1)}} \left\| \prod_{j=1}^{n+1} |e^{it\Delta}f_j|^{\frac{1}{n+1}}  \right\|_{L^{\frac{2(n+1)}{n}}(B_{2R})} \\
	& \lesssim K^{O(1)} M^{-\frac{1}{2(n+1)}} R^\epsilon\|f\|_{L^2}  \\
	& \leq K^{O(1)} M^{-\frac{1}{n+1}} \nu^{\frac{1}{n+1}} R^{\frac{\lambda}{2(n+1)}+\epsilon} \|f\|_{L^2} .
\end{align*}
Here in the fifth line we apply the change of variables  $x+R^{1-2\alpha}t^\alpha\mu(t/R^2)\rightarrow x$, and the last line holds due to $M\leq K^{O(1)} \nu^2 R^\lambda$ which can be derived by the definition of $\nu$.

\vskip0.5cm

\noindent \textbf{Narrow case.} We decompose $Y_{\text{narrow}}=Y_1 \cup Y_2$, where
$$  Y_1=Y_{\text{narrow}} \cap (B^n_R \times [0,K^{\frac{2}{1-\alpha}}]),\quad \quad Y_2=Y_{\text{narrow}} \backslash Y_1.    $$
Then
$$   	\|T_\alpha^{\gamma_R} f \|_{L^p(Y_{\text{narrrow}})} \leq  	\|T_\alpha^{\gamma_R} f \|_{L^p(Y_{1})} +	\|T_\alpha^{\gamma_R} f \|_{L^p(Y_{2})}.   $$
If the contribution from $Y_1$ dominates, this implies $Y_1$ contains $\sim ~M$ many $K^2$-cubes $B$. On the other hand, by the stationary phase method, $T_\alpha^{\gamma_R} f $ restricted on $Y_1$ is essentially supported on $B(0,K^{\frac{3}{1-\alpha}})\times [0,K^{\frac{2}{1-\alpha}}]$. Therefore, up to a rapidly decaying term, we have
 $$    \|T_\alpha^{\gamma_R} f \|_{L^2(Y_{1})} \lesssim K^{O(1)} \|f\|_{L^2}\lesssim K^{O(1)} \nu^{\frac{1}{n+1}} R^{\frac{\lambda}{2(n+1)}} \|f\|_{L^2} ,   $$
due to $\nu \geq K^{-2\lambda}$ and $\lambda \geq 1$. And then
 $$    \|T_\alpha^{\gamma_R} f \|_{L^p(Y_{1})} \lesssim K^{O(1)}M^{-\frac{1}{n+1}} \nu^{\frac{1}{n+1}} R^{\frac{\lambda}{2(n+1)}} \|f\|_{L^2} . $$

 Next we assume that the contribution from $Y_2$ dominates, and  this implies $Y_2$ contains $\sim ~M$ many $K^2$-cubes $B$. We break $B^n(0,R)$ in the physical space into $R/K$-cubes $D$. We write $f=\sum_{\tau,D} f_{\tau,D}$, where each $f_{\tau,D}$ is essentially supported on $D$ and Fourier supported on $\tau$. Then  $T_\alpha^{\gamma_R} f=\sum_{\tau,D} T_\alpha^{\gamma_R} f_{\tau,D}.$ Through the stationary phase method again, we get each $T_\alpha^{\gamma_R} f_{\tau,D}$ restricted on $B^\ast_R$ is essentially supported on
 $$     \Box:=\Box_{\tau,D}=\{ 0<t \leq R: |x-c_D+2tc_\tau|\leq R/K   \}, $$
where $c_D$ and $c_\tau$ denote the center of $D$ and $\tau$, respectively. We rewrite
$$|T_\alpha^{\gamma_R}f|=\left|    \int e\left[  x\cdot\xi+t\left|\xi+\frac{1}{2} R^{1-2\alpha}t^{\alpha-1} \mu\left(\frac{t}{R^2}\right)\right|^2 \right]  \widehat{f}(\xi) d\xi \right|.$$
We see that on $Y_2$,
$$   \left(\xi,\left|\xi+\frac{1}{2} R^{1-2\alpha}t^{\alpha-1} \mu\left(\frac{t}{R^2}\right)\right|^2\right) \subset N_{K^{-2}}(\mathbb{P}^n).  $$
 Therefore, using the lower dimensional $l^2$ decoupling inequality of variable coefficient version  \cite{ILX}, we get
\begin{equation}\label{dec}
\|T_\alpha^{\gamma_R}f\|_{L^p(B)} \lesssim K^{\epsilon^2} \left(  \sum_{\tau \in \mathcal{S}(B)}\|T_\alpha^{\gamma_R}f_\Box\|^2_{L^p(w_B)}  \right)^{\frac{1}{2}},
\end{equation}
where $w_B$ is given by
$$   w_B(x):=  \left(  1+\frac{|x-c_B|}{K^2} \right)^{-100n},   $$
and $c_B$ denotes the center of $B$.

Next we will perform several dyadic pigeonholing which has similar spirit as in \cite[Section 3]{DZ}. Set $R_1=R/K^2=R^{1-2\delta}$, $K_1=R_1^\delta=R^{\delta-2\delta^2}$. Tile $\Box$ by $KK_1^2 \times ... \times KK_1^2 \times K^2K_1^2$-tubes $S$ running parallel to the long axis of $\Box$. We now perform dyadic pigeonholing argument to $S$ and $\Box$:

\noindent (1) For each $\Box$, we sort $S \subset \Box$ that intersect $Y_2$ according to the value $\|T_\alpha^{\gamma_R} f\|_{L^p(S)}$ and the number of $K^2$-cubes in $Y_2$ contained in it: for dyadic numbers $\eta,\beta_1$, define
	$$    \mathbb{S}_{\Box,\eta,\beta_1}:=\{ S\subset \Box:S~\text{contains}\sim\eta~\text{many}~K^2\text{-cubes~in~}Y_{2},~ \|T_\alpha^{\gamma_R} f\|_{L^p(S)}\sim\beta_1  \}.    $$
Let $Y_{\Box,\eta,\beta_1}$ be the union of the tubes $S$ in $\mathbb{S}_{\Box,\eta,\beta_1}$.  	
	
\noindent (2) For fixed $\eta$ and $\beta_1$, we sort $\Box$ according to the value $\|f_\Box\|_{L^2}$, the number $\#\mathbb{S}_{\Box,\eta,\beta_1}$ and the value $\nu_1$ given below: for dyadic numbers $\beta_2,M_1,\nu_1$, define
	$$   \mathbb{B}_{\eta,\beta_1,\beta_2,M_1,\gamma_1}:=\left\{ \Box: \|f_\Box\|_{L^2}\sim \beta_2,~ \# \mathbb{S}_{\Box,\eta,\beta_1}\sim M_1,~ \max_{T_r\subset \Box:r\geq K_1^2}    \frac{\#\{ S\in \mathbb{S}_{\Box,\eta,\beta_1}:S \subset T_r  \}}{r^\lambda}      \sim\nu_1 \right\},  $$
	where $T_r$ are $Kr \times ...\times Kr \times K^2 r$-tubes in $\Box$ running parallel to the long axis of $\Box.$

On each $B \subset Y_2$, up to a rapidly decaying term, we have
$$ T_\alpha^{\gamma_R} f =\sum_{\eta,\beta_1,\beta_2,M_1,\nu_1}\left( \sum_{\substack{\Box\in \mathbb{B}_{\eta,\beta_1,\beta_2,M_1,\nu_1}  \\ B\subset Y_{\Box,\eta,\beta_1}}} T_\alpha^{\gamma_R} f_\Box  \right). $$
Under the hypothesis $\|f\|_{L^2}=1$, we can further assume that
\begin{align*}
	1\leq \eta \leq& K^{O(1)}, \quad R^{-C} \leq \beta_1 \leq K^{O(1)}, \quad R^{-C} \leq \beta_2 \leq 1,    \\
	& 1 \leq M_1 \leq R^{O(1)},\quad K^{-2n} \leq \nu_1 \leq R^{O(1)}
\end{align*}
for some constant $C$. Thus there are only $O(\log R)$ choices for each dyadic number. By dyadic pigeonholing, these exist $\eta,\beta_1,\beta_2,M_1,\nu_1$ depending on $B$ such that
\begin{equation}\label{absd}
	\|T_\alpha^{\gamma_R} f\|_{L^p(B)} \lesssim (\log R)^5 \left\|  \sum_{\Box \in\mathbb{B}_{\eta,\beta_1,\beta_2,M_1,\nu_1} }T_\alpha^{\gamma_R} f_\Box \cdot \chi_{Y_{\Box,\eta,\beta_1}}  \right\|_{L^p(B)}.
\end{equation}

Finally, we sort $B \subset Y_2$. Since there are only $O(\log R)$ choices on $\eta,\beta_1,\beta_2,M_1,\nu_1$, by pigeonholing, we can find uniform  $\eta,\beta_1,\beta_2,M_1,\nu_1$ such that (\ref{absd}) holds for a fraction $\gtrsim (\log R)^{-6}$ of all $B \subset Y_2$. We denote the union of such $B$ by $Y'_2$. From now on, we abbreviate $Y_{\Box,\eta,\beta_1}$ and $\mathbb{B}_{\eta,\beta_1,\beta_2,M_1,\nu_1}$ to $Y_\Box$ and $\mathbb{B}$, respectively. Next we further sort $B \subset Y_2'$ by the number $\# \{ \Box \in \mathbb{B}: B \subset Y_2'  \}$: for dyadic number $\ell$, define
$$  \mathbb{Y}_\ell :=\{  B \subset Y'_{2}:  \#  \{ \Box \in \mathbb{B}: B\subset Y_\Box \}\sim\ell  \}.  $$
Let $Y_{\ell}$ be the union of $B$ in $\mathbb{Y}_\ell$. Using dyadic pigeonholing again, we can choose $\ell$ such that
\begin{equation}\label{flr eq7}
	\#\{ B : B \subset Y_{\ell} \} \gtrsim (\log R)^{-6} \#\{B : B \subset Y'_{2} \}.
\end{equation}
From now on, we fix $\ell$, and denote $Y_{\ell}$ by $Y''$.

Now on each $B \subset Y''$, by (\ref{dec}), we obtain
\begin{align*}
	\| T_\alpha^{\gamma_R} f \|_{L^p(B)} \lesssim&  (\log R)^5 \left\|  \sum_{\Box \in\mathbb{B} }T_\alpha^{\gamma_R} f_\Box \cdot \chi_{Y_{\Box}}  \right\|_{L^p(w_B)}  \\
	\lesssim  & (\log R)^5  K^{\epsilon^2} \left( \sum_{
		\Box \in \mathbb{B}: B \subset Y_\Box} \|T_\alpha^{\gamma_R} f_\Box \|^2_{L^p(w_B)}  \right)^{\frac{1}{2}}  \\
	\lesssim   &  (\log R)^5 K^{\epsilon^2}\ell^{\frac{1}{n+1}} \left( \sum_{\Box \in \mathbb{B}: B \subset Y_\Box} \|T_\alpha^{\gamma_R} f_\Box \|^p_{L^p(w_B)}  \right)^{\frac{1}{p}}.
\end{align*}
Summing $B \subset Y''$, note (\ref{flr eq7}), we see that
\begin{align}
	\|T_\alpha^{\gamma_R} f \|_{L^p(Y)} &\lesssim (\log R)^6 \|T_\alpha^{\gamma_R} f \|_{L^p(Y'')}     \nonumber \\
	& \lesssim  (\log R)^{11} K^{\epsilon^4} \ell^{\frac{1}{n+1}} \left(\sum_{\Box \in \mathbb{B}} \|T_\alpha^{\gamma_R} f _\Box\|^p_{L^p(Y_\Box)}\right)^{\frac{1}{p}} .   \label{11}
\end{align}
As for the term $ \|T_\alpha^{\gamma_R} f _\Box\|_{L^p(Y_\Box)}$, we will use
parabolic rescaling and induction on scale. For each $\tau=\tau_\Box$, write $\xi=c_\tau+K^{-1}\eta$, then
$$   |T_\alpha^{\gamma_R} f _\Box(x,t)|=K^{-\frac{n}{2}}|T_\alpha^{\gamma_{R/K}} g(\tilde{x},\tilde{t})   | $$
for some $g$ with supp$\widehat{g}\subset B^n(0,1)$ and $\|g\|_{L^2}=\|f_\Box\|_{L^2}$, and $\tilde{x},\tilde{t}$ are given by $\tilde{x}=K^{-1}x+2tK^{-1}c_\tau$, $\tilde{t}=K^{-2}t$. Thus
\begin{equation}\label{12}
 \|T_\alpha^{\gamma_R} f _\Box(x,t)\|_{L^p(Y_\Box)}=K^{-\frac{1}{n+1}} \|T_\alpha^{\gamma_{R/K}}g(\tilde{x},\tilde{t})\|_{L^p(\tilde{Y})},
\end{equation}
where $\tilde{Y}$ is the image of $Y_\Box$ under the new coordinates. We can see that $\|T_\alpha^{\gamma_{R/K}}g\|_{L^p(\tilde{Y})}$ just satisfies the condition of Proposition \ref{pro} at scale $R_1$ under the new parameters $M_1, \nu_1, R_1$. Therefore we apply inductive hypothesis to obtain
\begin{equation}\label{13}
\|T_\alpha^{\gamma_{R/K}}g\|_{L^p(\tilde{Y})} \lesssim  M_1^{-\frac{1}{n+1}} \nu_1^{\frac{1}{n+1}} R_1^{\frac{\lambda}{2(n+1)}+\epsilon} \|g\|_{L^2}.
\end{equation}
Combining (\ref{11}), (\ref{12}), (\ref{13}), one has
\begin{align}
	\|T_\alpha^{\gamma_R} f \|_{L^p(Y)} & \leq  K^{2\epsilon^4}  \ell^{\frac{1}{n+1}} K^{-\frac{1}{n+1}} M_1^{-\frac{1}{n+1}} \nu_1^{\frac{1}{n+1}} R_1^{\frac{\lambda}{2(n+1)}+\epsilon} (\sum_{\Box \in \mathbb{B}}  \|f_\Box\|_{L^2}^p )^{\frac{1}{p}}\nonumber   \\
	& \lesssim   K^{2\epsilon^4}  (\frac{\ell}{\#\mathbb{B}})^{\frac{1}{n+1}} K^{-\frac{1}{n+1}} M_1^{-\frac{1}{n+1}} \nu_1^{\frac{1}{n+1}} R_1^{\frac{\lambda}{2(n+1)}+\epsilon}  \|f\|_{L^2}. \label{14}
\end{align}

Finally, we have the relations on the old parameters $M,\nu$ and new parameters $M_1,\nu_1$ as follows:
\begin{equation}\label{15}
\frac{\ell}{\#\mathbb{B}}\lesssim \frac{(\log R)^6M_1\eta}{M},\quad \quad\eta\lesssim\frac{\nu K^{\lambda+1}}{\nu_1}.
\end{equation}
(\ref{15}) can be proved through the same steps as in \cite[Section 3]{DZ}, here we omit relevant details. We add (\ref{15}) to the estimate (\ref{14}), then
$$    \| T_\alpha^{\gamma_R} f \|_{L^p(Y)}  \lesssim K^{2\epsilon^2-2\epsilon}M^{-\frac{1}{n+1}} \nu^{\frac{1}{n+1}}R^{\frac{\lambda}{2(n+1)}+\epsilon} \|f\|_{L^2}.$$
Since $K=R^\delta$, we can choose sufficiently large $R$, then the induction closes.

\qed

Finally, we say a few words on the case $0 <\alpha <1/2$. In Proposition \ref{pro}, the condition $1/2 \leq \alpha<1$ is used several times, such as locally constant property in the broad case, and the $l^2$ decoupling inequality in the narrow case. These features seem to show that $\alpha=1/2$ is the critical point, and the pointwise convergence will become worse when $\alpha<1/2$. In fact, if $0 <\alpha <1/2$, we still can reduce Theorem \ref{th1} to Theorem \ref{th4} via the same argument. But the relevant properties of $T_\alpha^{\gamma_R} f$ become worse. For example, if $f$ is Fourier supported on $B(0,1)$, then $|T_\alpha^{\gamma_R} f|$ can be viewed as constant on the box at scale $1\times ...\times 1 \times R^{2-1/\alpha}$. On the other hand, if we consider its wave packet decomposition,  write $T_\alpha^{\gamma_R}f=\sum_{\theta,\nu}T_\alpha^{\gamma_R}f_{\theta,\nu}$, where $\theta$ denotes balls of radius $R^{-1/2}$ in the frequency space, $\nu$ denotes balls of radius $R^{1/2}$ in the physical space, then each $T_\alpha^{\gamma_R}f_{\theta,\nu}$ restricted on $B^\ast_R$ is essentially supported on the tube at scale $R^{1-\alpha} \times ...\times R^{1-\alpha} \times R$. These features make the same argument on Proposition \ref{pro} fail. However, for $n=1$, Cho-Lee-Vargas \cite{CLV} proved (\ref{ntc}) holds for $s>\max(1/2-\alpha,1/4)$ if $0<\alpha < 1$. In other words, the critical point when $n=1$ is $\alpha=1/4$. Therefore, it seems that new ingredients are still needed for the case $0<\alpha<1/2$ when $n\geq 2$.

\subsection*{Acknowledgements}   This project were supported by the National Key R\&D Program of China: No. 2022YFA1005700 and 2020YFAO712903,  and NSFC No. 11831004.

\bibliographystyle{plain}
\bibliography{mybibfile}

\end{document}